\documentclass[conference]{IEEEtran}
\IEEEoverridecommandlockouts
\usepackage{cite}
\usepackage{amsmath,amssymb,amsfonts}
\usepackage{algorithmic}
\usepackage{graphicx}
\usepackage{hyperref}
\usepackage{textcomp}
\usepackage{xcolor}
\usepackage{verbatim}
\usepackage{float}
\usepackage{longtable,supertabular,subcaption}
\usepackage{pgfplots}
\usepackage[export]{adjustbox}
\usepackage{enumitem}
\usepackage{adjustbox}

\def\BibTeX{{\rm B\kern-.05em{\sc i\kern-.025em b}\kern-.08em
    T\kern-.1667em\lower.7ex\hbox{E}\kern-.125emX}}
\begin{document}

\title{Land Use Detection \& Identification using Geo-tagged Tweets\\
}

\author{\IEEEauthorblockN{Saeed Khan}
\IEEEauthorblockA{\textit{School of ITEE} \\
\textit{University of Queensland}\\
Brisbane, Australia \\
s.khan@uq.edu.au}
\and
\IEEEauthorblockN{Md Shahzamal}
\IEEEauthorblockA{\textit{Department of Computing} \\
\textit{Macquarie University}\\
Sydney, Australia \\
md.shahzamal@mq.edu.au}


}

\maketitle

\begin{abstract}
Geo-tagged tweets can potentially help with `sensing' the interaction of people with their surrounding environment. Based on this hypothesis, this paper makes use of geo-tagged tweets in order to ascertain various land uses with a broader goal to help with urban/city planning. The proposed method utilises supervised learning to reveal spatial land use within cities with the help of Twitter activity signatures. Specifically, the technique involves using tweets from three cities of Australia namely Brisbane, Melbourne and Sydney. Analytical results are checked against the zoning data provided by respective city councils and a good match is observed between the predicted land use and existing land zoning by the city councils. We show that geo-tagged tweets contain features that can be useful for land use identification. \newline
\end{abstract}


\section{Introduction}
The goal of urban planning is to develop cities with the best possible productive use of land, the ultimate purpose being improving peoples' day to day life and better management of resources. Generally, the urban planning characterises residential, commercial or recreational use of the land. The land classification or zoning determines how each segment of the land is to be used. Zoning can be defined as the process of dividing land in a city council (municipality) into various segments based upon its intended use. Thus, it may enforce a variety of specific or conditional uses of land. It may also indicate the size and dimensions of land area as well as the form and scale of constructions. One of the issues associated with zoning is that it requires actual assessment with regards to the use of land as originally planned. This information needs to be collected on site, and is usually obtained through survey or questionnaire to record how people interact with their land. This approach has restrictions in the sense that it is very much dependent upon the willingness of people to participate in such activity, involves costs, and, essentially is just a snapshot in time. Technologies such as GIS \cite{yin2011monitoring} may reveal some of the uses of land through image processing. However, since images are not acquired very frequently, such techniques may not be very productive. \newline

As a result of using a large number of mobile devices, people leave behind digital footprints of their interaction with the environment. Availability of a vast number of mobile applications has turned the mobile devices into sensors of human activity in space and time. This has led to the advent of new research areas such as smart cities and urban computing with a goal to understand the city dynamics through ubiquitous data, and thus improving the quality of life. Data sources such as GPS, Bluetooth, Wifi and cell phone records are becoming more and more useful for applications such as transport planning \cite{Frias-Martinez:2012:EUC:2346496.2346499}, traffic estimation \cite{6172682} and social studies \cite{oloritun2013identifying}. This paper suggests a method to use Twitter data to determine the type of land people tweet from by analysing the temporal patterns of their tweeting behaviour. The method is scalable and efficient compared to some other techniques in the literature and can utilise any number of data points in order to capture the best possible spatial features.

\section{Related Work}
For urban planning, a number of different data sources have been used by various researchers in their studies to model and characterise land use. The data for these studies comes from either a GPS device, cell phone, or a Location Based Social Network service. A quick review of work related to each of these technologies is provided here.\newline

\paragraph{GPS:} This data source has limited availability, is recorded every few seconds and is usually based on a small number of subjects. It is not very widely used for land classification purposes. The data set is usually acquired through buses and taxis e.g.\ in \cite{yuan2012discovering}, the authors conduct analysis of GPS tracks obtained from taxis to infer the mobility of individuals and detect areas of various activities using a pre-defined map of the places of interest. \newline

\paragraph{Cell Phone Data:} Usually known as Call Detail Records (CDRs) contain information for a large number of people, but is not easy to acquire because of privacy issues and with location precision at the level of mobile tower only when some activity is initiated. In \cite{5594641}, the authors use cell phone data to analyse the relationship between cellular activity and commercial land use to identify the prevailing pattern in each area. They qualitatively presented their results and did not use any land use information. In another similar study \cite{soto2011robust}, the authors analysed the cell phone activity to characterise and cluster similar locations. \newline

\paragraph{Location Based Social Networks:} These offer precise location information in the form of latitude and longitude when users enable the geo-location feature. This data set can be most up to date and is easy to acquire as a result of ubiquitous use of various mobile and social network applications by people. \cite{noulas2011exploiting} model activity patterns in New York City and London by using the geo-located information provided by a social networking platform. They characterise the activity patterns by using a set of pre-defined categories which indicate the type of check-in locations. Due to this approach, their results provide a coarse idea of the land use. In another work, \cite{cranshaw2012livehoods} present a method to explore land use and user activities at a large scale with the help of Foursquare platform. They validate the results with personal interviews to confirm the identified land uses. In \cite{Garcia}, the authors group twitter users by city zone and various times of the day to study the daily dynamics of the city. For a zone, they compare daytime activity with night time to determine zones showing increased activity. Using multiple regression, they determine how various land uses influence different time of the day via changing coefficients. They observed that user activity seemed to decrease throughout the day for land uses such as education, health and offices, remained constant for parks, and increased for retail and residential zones. This study aims to establish the link between twitter activity and the known land use types and does not intend to determine the land uses from twitter activity itself. In \cite{Lee2012U} and \cite{Wakamiya}, the authors propose a method to characterise geographic regions by focusing on day-to-day living patterns of crowds using Twitter data. To grasp images of a city, the
regularity of a region is defined by using activity patterns such as the number of tweets,
their contents, the number of users and the movement of crowds. Moreover, the changes
in regularity patterns with respect to time are analysed and the grouped urban types are
categorised by tracking similar patterns across the regions.
\newline

In this study, we also make use of publicly available geo-tagged tweets for the identification of land use in a city. However, we purely use the temporal patterns of raw twitter activity without using any textual information or user profiles. Our method is simplistic, requires less processing and thus can be applied on any volume of data. Specifically, we use temporal signatures of Brisbane tweets as a template (training data) against which land use in other cities is subsequently detected and predicted. We validate the results against land use data provided by councils of test cities and observe a good match between the predicted land use and actual land classification. \newline

\section{Data Set \& Experimental Setup}
The data set has been collected using Twitter Streaming API \cite{Twitter} covering Australia and spanning a period from January 2015 till July 2016. It consists of 9462345 tweets from 245796 unique users.
For Sydney we have 947964 geo-tagged tweets and 40281 unique users, for Melbourne we have 854821 geo-tagged tweets and 35556 unique users and for Brisbane we have 276394 geo-tagged tweets and 14555 unique users in the data set. The proposed method involves using the spatial location of tweets for above three cities of Australia. Brisbane is chosen as the baseline city, and Melbourne and Sydney are used as test cities on which the method is applied. For Brisbane, Central Business District (CBD) is chosen as the reference around which various land use zones are identified based on its city council's scheme.\footnote{Schemes can be found \href{https://www.brisbane.qld.gov.au/planning-and-building/planning-guidelines-and-tools/brisbane-city-plan-2014/brisbane-city-plan-2014-mapping/zoning-maps}{\color{blue}here}.}  Although the council divides land into numerous types, key zones consist of business, residential, education and recreation since these zones account for the highest twitter activity compared to other zones. \newline

Starting from the Brisbane CBD, corresponding temporal signatures are extracted for each key zone. These signatures are later used as a template. The hypothesis is that areas in different cities can be classified by their similarity to these temporal signatures. To detect land use in our test cities, CBDs are selected again and various activity clusters are generated from the city-blocks until some meaningful temporal patterns are obtained. Then, the closest matching temporal pattern from Brisbane is found for each unknown cluster and the land use is thus predicted. This predicted cluster is then compared to actual land use based on its city zoning scheme and the overlap percentage is calculated for evaluation. \newline

\section{Temporal Signature Generation}
\begin{figure*}[htbp]
\centering
\includegraphics[width=17cm]{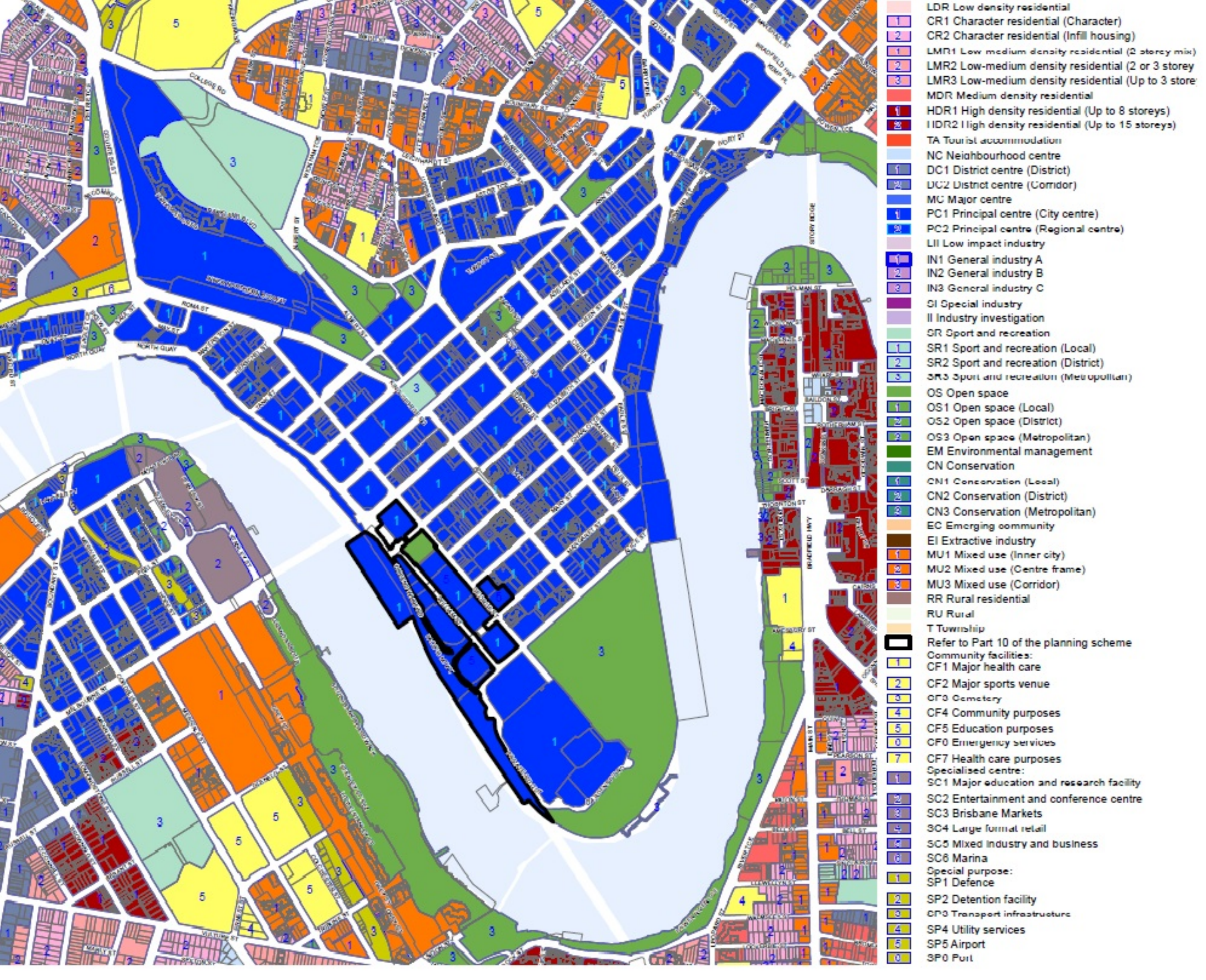}
\caption{Brisbane City Zoning Map: land use is identified by a combination of colour and number e.g.\ business (Blue 1), residential (Light pink 1, 2, 3 \& Maroon 1, 2), education (Yellow 5) and recreation (Yellow 2 \& Grey 2). \newline
Source: cityplan.brisbane.qld.gov.au/eplan/\#/Property/0}
\label{fig:brisbanecityzoningmap}
\end{figure*}

Using Brisbane as baseline city,
each of the key identified zone is analysed for its twitter activity. Figure \ref{fig:brisbanecityzoningmap} shows the city council's land zoning scheme based on areas surrounding the CBD.
Since each geo-tagged tweet contains its location information and time-stamp, the temporal pattern for every hour of the day is extracted and plotted. This temporal signature reflects the twitter activity for a given land zone and spans the entire period of the data set. As a result, four temporal signatures are generated for business (Blue 1), residential (Light pink 1, 2, 3 \& Maroon 1, 2), education (Yellow 5) and recreation (Yellow 2 \& Grey 2) land use, respectively. \newline

\begin{figure}[htbp]
\centering
\captionsetup{format= hang, justification=centering}
  \begin{subfigure}[b]{0.53\textwidth}
         \includegraphics[width=\textwidth]{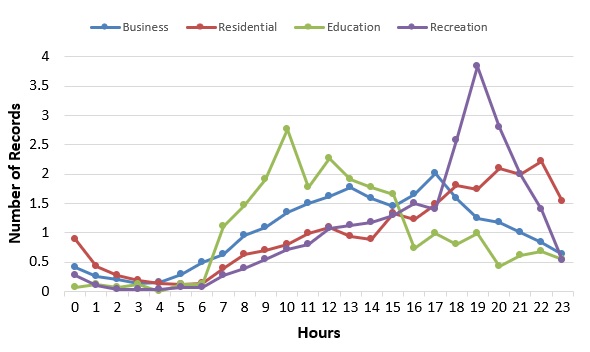}
  \end{subfigure}

\caption{Brisbane Temporal Signatures}
\label{fig:brisbanetempsigs}
\end{figure}

Figure \ref{fig:brisbanetempsigs} shows the temporal behaviour of each key land zone as the normalised hourly tweets for the entire data set. Normalisation is the process of adjusting values from different scales to some common scale in order to provide a better comparison of the values. For example, the raw data for Business cluster consists of very high values for each hour compared to other clusters, and gives an impression that it is probably the most `active' cluster. However, a better analysis would be to figure out the average percentage of tweets posted in each hour of a day within each of the clusters. Consequently, each hourly record value is divided by the overall average of all hours to get the corresponding normalised value. In this way, the average for all normalised values is 1 for each cluster.\newline
\begin{equation*}
    T_h^c = \frac{\#T_h^c}{\sum_{h=0}^{23} \#T_h^c} 
\end{equation*}
where $T_h^c$ is the percentage of tweets in hour $h \in \{0,\dots,23\}$ within cluster $c$, and $\#T_h^c$ is the raw count of such tweets.

As shown in Figure \ref{fig:brisbanetempsigs}, the business zone is characterised by a large tweeting activity during day time. The first peak in the tweeting activity is observed around 1 p.m. and the second peak is seen at 5 p.m. This might be related to the times corresponding to activities such as going for a lunch or going home. Spatially, this zone covers areas such as hotels, restaurants, shopping centres and offices. Some bus and train stations are also included in it. Hence, overall this zone represents `business' activity in the city. The residential zone shows increasing twitter activity in the afternoon with maximum shown at 10 p.m. It covers areas comprising of private residences, apartments and hostels. The education zone has its maximum peak at 10 a.m. followed by another one at 12 p.m. After this, the trend in the twitter activity gradually declines. This zone covers area including a school and a TAFE in South Brisbane. The recreation zone has its distinctive peak at 7 p.m. after which there is a sharp decline in twitter activity. Spatially it covers areas including Suncorp Stadium and Brisbane Convention \& Exhibition Centre. We hypothesise that Twitter activities in other cities such as Melbourne and Sydney have the same pattern for each of the four land uses. We test this hypothesis through our proposed land use detection model. 

\section{Land Use Detection} 
In this section, land use is detected for each test city: Melbourne and Sydney by using the spatial location of tweets.
The cluster building process is manual and incremental as per following steps: 
\begin{enumerate}[wide, labelwidth=!, labelindent=0pt]
\item Start from CBD, randomly pickup a city-block, get its latitude-longitude bounding box, retrieve the data points contained inside and plot its temporal curve. If the curve contains data points for every hour of the day, keep the city-block as one of the clusters. Otherwise:

\item Manually expand the cluster by adding adjoining areas by enlarging the latitude-longitude bounding box. Check the temporal pattern again and expand further if necessary until a full temporal curve showing data for every hour of the day is obtained.

\item Choose another city-block and generate clusters in above manner.

\item Discard any cluster with incomplete temporal curve i.e missing data points for certain hour(s) of the day.

\end{enumerate}

This cluster generation process emphasises that there should be no zero reading for any hour because this might indicate that such a cluster is either just random or irregular in nature with not much activity taking place in it. In comparison, a cluster which has activity for every hour of the day shows that it is consistently and regularly used by people and this is important for determining an activity or land use eventually. It should also be noted that some of the clusters may eventually form a little far from the CBD as a result of this process.

\subsection{Melbourne Clusters}
\begin{figure*}[!h]
\centering
  \begin{subfigure}[b]{0.38\textwidth}
      \centering
        \includegraphics[width=6.0cm,height=5.1cm]{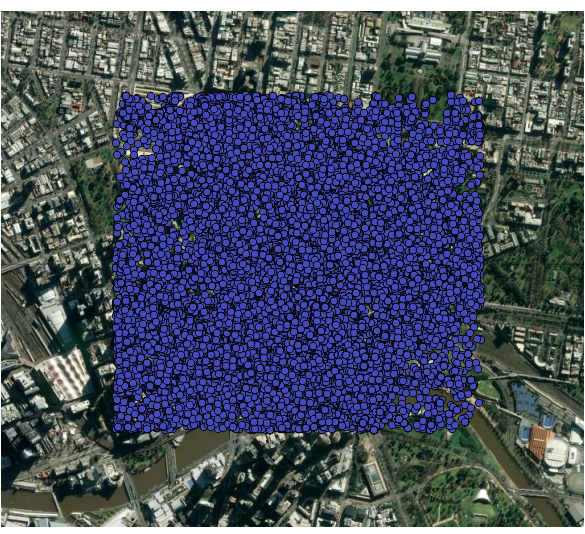}
         \caption{Cluster 1}
  \end{subfigure} 
  \begin{subfigure}[b]{0.38\textwidth}
      \centering
        \includegraphics[width=6.0cm,height=5.1cm]{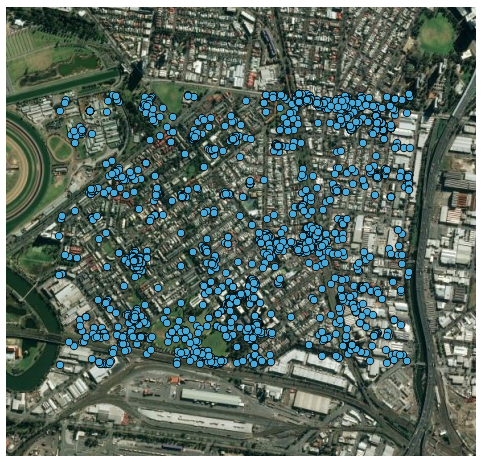}
        \caption{Cluster 2}
  \end{subfigure} \\ \vspace{0.4cm}
  \begin{subfigure}[b]{0.38\textwidth}
      \centering
         \includegraphics[width=6.0cm,height=5.09cm]{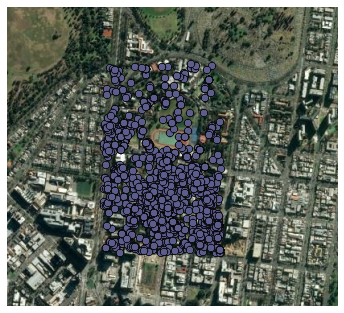}
         \caption{Cluster 3}
  \end{subfigure} 
  \begin{subfigure}[b]{0.38\textwidth}
      \centering
        \includegraphics[width=6.0cm,height=5.1cm]{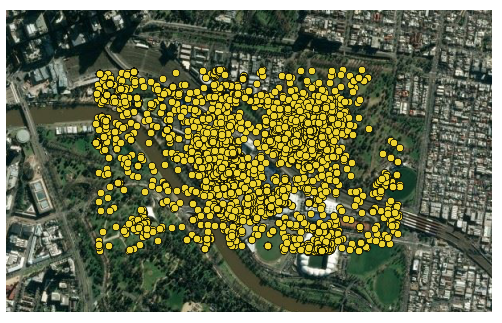}
        \caption{Cluster 4}
  \end{subfigure} 
\caption{Melbourne Clusters. Cluster 1 ranges from -37.820565 to -37.806470 (latitude) and 144.955165 to 144.974617 (longitude). Cluster 2: -37.799383 to -37.787955 (latitude) and 144.916789 to 144.935543 (longitude). Cluster 3: -37.800277 to -37.792296 (latitude) and 144.957810 to 144.964049 (longitude). Cluster 4: -37.824938 to -37.817277 (latitude) and 144.971911 to 144.988133 (longitude).}
\label{fig:melbclusters}
\end{figure*}

For Melbourne, the cluster building process explained above is followed and eventually four clusters are generated which are then subjected to further processing and analysis. These clusters are shown in figure \ref{fig:melbclusters}. 

\begin{figure}[h!]
\centering
\captionsetup{format= hang, justification=centering}
  \begin{subfigure}[b]{0.52\textwidth}
         \includegraphics[width=\textwidth]{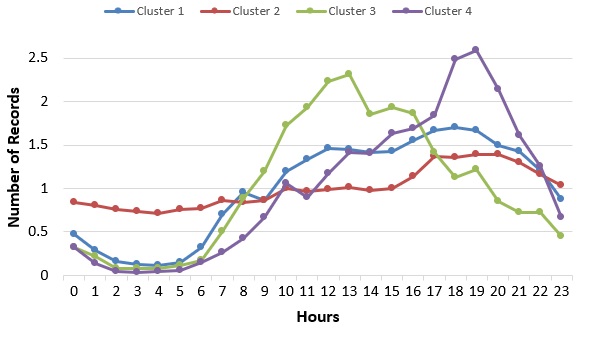}
  \end{subfigure}

\caption{Melbourne Temporal Signatures}
\label{fig:melbtempsigs}
\end{figure}

The temporal patterns for each of these unknown clusters are shown in figure  \ref{fig:melbtempsigs} and these plots represent normalised data. 
Each of these clusters is yet unknown in terms of the land use that it actually represents. As a next step, these clusters are compared with Brisbane's known zones to judge their similarity. This similarity is evaluated quantitatively by using the measure of mean-squared error. The mean-squared error (MSE) function calculates the average of squares of errors, or in other words, the average squared difference between the observed values (here a temporal activity signature from Brisbane) and predicted values (here a temporal activity signature from Melbourne). It is always non-negative and values closer to zero are better as they indicate less error, with zero signifying no error at all. 
\newcommand{\MSE}{\ensuremath{\,\textrm{MSE}}}
\begin{equation}
\MSE =  \ \frac {1}{n} \sum_{i=1}^{n}(X_i - Y_i)^2
\end{equation}

\noindent where \textit{n} is the number of data points,$\ \textit{X}_i$ represents observed values and$\ \textit{Y}_i$ represents predicted values. The similarity between values can also be measured with other techniques such as Principal Component Analysis (PCA) \cite{pearson1901liii}. PCA is based on dimensionality reduction concept and involves building vectors based on observed and predicted values. This work, however, uses MSE technique because of its simplicity in calculating errors. Each extracted cluster of Melbourne is checked against every zone of Brisbane and the corresponding mean-squared error is calculated. This error is calculated on the basis of each hourly normalised value for a cluster and then the overall error is summed. Table \ref{tab:melbmsetable} lists errors for each pair of Melbourne and Brisbane cluster below: 

\begin{table}[h]
  \centering
    \begin{tabular}{ccccc} \hline
    {\textbf{}} & \multicolumn{4}{c}{\textbf{Brisbane Zones}} \\
\multicolumn{1}{c}{} & \multicolumn{1}{c}{\textbf{Business}} & \multicolumn{1}{c}{\textbf{Residential}} & \multicolumn{1}{c}{\textbf{Education}} & \multicolumn{1}{c}{\textbf{Recreation}} \\ \hline
    
 Melb. Cluster 1 & \textbf{0.045}  & 0.156  & 0.419  & 0.387 \\ 
    
 Melb. Cluster 2 & 0.206  & \textbf{0.204}  & 0.611  & 0.607 \\
    
 Melb. Cluster 3 & \textbf{0.098}  & 0.613  & 0.184  & 0.90 \\
    
 Melb. Cluster 4 & 0.239  & 0.190  & 0.834  & \textbf{0.117} \\\hline
    \end{tabular}
\caption{MSE between Melbourne \& Brisbane reference clusters}
\label{tab:melbmsetable}
\end{table}

The minimum error for each pair of clusters is highlighted in bold and based on this, the cluster for Melbourne is assigned. e.g.\ Cluster 1 has the least error when compared against the Business zone of Brisbane. So the model has predicted that Cluster 1 of Melbourne probably belongs to Business area of the city. Similarly, Cluster 2 of Melbourne has least error against Residential zone of Brisbane and so this cluster is predicted as Residential. In this way, a land use prediction is made for all Melbourne clusters to be later validated against actual land use based on its city council's classification scheme. If the predictions are mostly correct, it is indicative that we can use a city's twitter activity pattern as a reference to identify the land use in other similar cities.

\subsection{Sydney Clusters} 
\begin{figure*}[htbp]
\centering
  \begin{subfigure}[b]{0.38\textwidth}
      \centering
         \includegraphics[width=6.0cm,height=5.1cm]{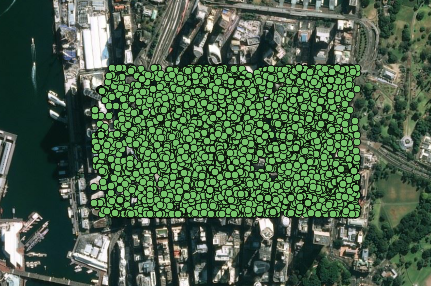}
         \caption{Cluster 1}
  \end{subfigure}
  \begin{subfigure}[b]{0.38\textwidth}
      \centering
        \includegraphics[width=6.0cm,height=5.1cm]{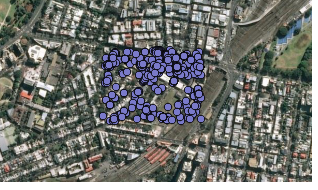}
        \caption{Cluster 2}
  \end{subfigure} \\ \vspace{0.4cm}
  \begin{subfigure}[b]{0.38\textwidth}
      \centering
         \includegraphics[width=6.0cm,height=5.1cm]{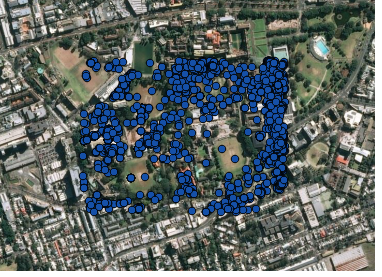}
         \caption{Cluster 3}
  \end{subfigure}
  \begin{subfigure}[b]{0.38\textwidth}
      \centering
        \includegraphics[width=6.25cm,height=5.22cm]{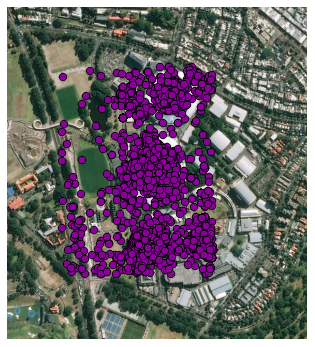}
        \caption{Cluster 4}
  \end{subfigure} 
    
\caption{Sydney Clusters. Cluster 1 ranges from -33.868316 to -33.863256 (latitude) and 151.202142 to 151.213085 (longitude). Cluster 2: -33.890915 to -33.888452 (latitude) and 151.196287 to 151.200373 (longitude). Cluster 3: -33.891788 to -33.886575 (latitude) and 151.182370 to 151.190743 (longitude). Cluster 4: -33.895203 to -33.888043 (latitude) and 151.220731 to 151.227083 (longitude).}
\label{fig:sydnclusters}
\end{figure*}

For Sydney, the clusters are built in a similar fashion and four clusters are generated for further processing and analysis. These clusters are shown in figure \ref{fig:sydnclusters}. 

\begin{figure}[htbp]
\centering
  \begin{subfigure}[b]{0.51\textwidth}
         \includegraphics[width=\textwidth]{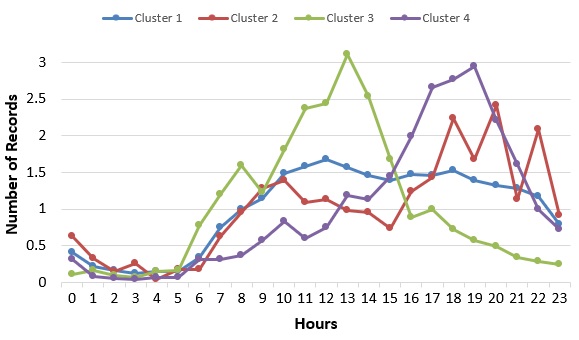}
  \end{subfigure}
  
\caption{Sydney Temporal Signatures}
\label{fig:sydntempsigs}
\end{figure}

The temporal patterns of each of these unknown clusters are shown in figure  \ref{fig:sydntempsigs}. The plots represent normalised hourly twitter activity for the corresponding cluster. 
These clusters are yet unknown in terms of the land use that they actually represent. Like Melbourne, these clusters are compared with Brisbane's known zones to judge their mutual similarity. Each Sydney cluster is matched against every zone of Brisbane and their corresponding mean-squared error is calculated. Table \ref{tab:sydnmsetable} lists errors for each pair of Sydney and Brisbane cluster. The minimum MSE error for each pair of clusters is highlighted and based on this, the cluster for Sydney is assigned. Cluster 1 is predicted as Business, Cluster 2 as Residential, Cluster 3 as Education, and Cluster 4 as Recreation. \newline

\begin{table}[htbp]
  \centering
    \begin{tabular}{ccccc} \hline
    {\textbf{}} & \multicolumn{4}{c}{\textbf{Brisbane Zones}} \\
\multicolumn{1}{c}{} & \multicolumn{1}{c}{\textbf{Business}} & \multicolumn{1}{c}{\textbf{Residential}} & \multicolumn{1}{c}{\textbf{Education}} & \multicolumn{1}{c}{\textbf{Recreation}} \\ \hline
    
 Sydney Cluster 1 & \textbf{0.030}  & 0.238  & 0.267  & 0.536 \\ 
    
 Sydney Cluster 2 & 0.266  & \textbf{0.114}  & 0.691  & 0.355 \\
    
 Sydney Cluster 3 & 0.398  & 1.204  & \textbf{0.198}  & 1.613 \\
    
 Sydney Cluster 4 & 0.397  & 0.321  & 1.190  & \textbf{0.149} \\\hline
    \end{tabular}
\caption{MSE between Sydney \& Brisbane reference clusters}
\label{tab:sydnmsetable}
\end{table}

\section{Land Use Validation}
\begin{figure*}[h!]
\centering
\includegraphics[width=17cm]{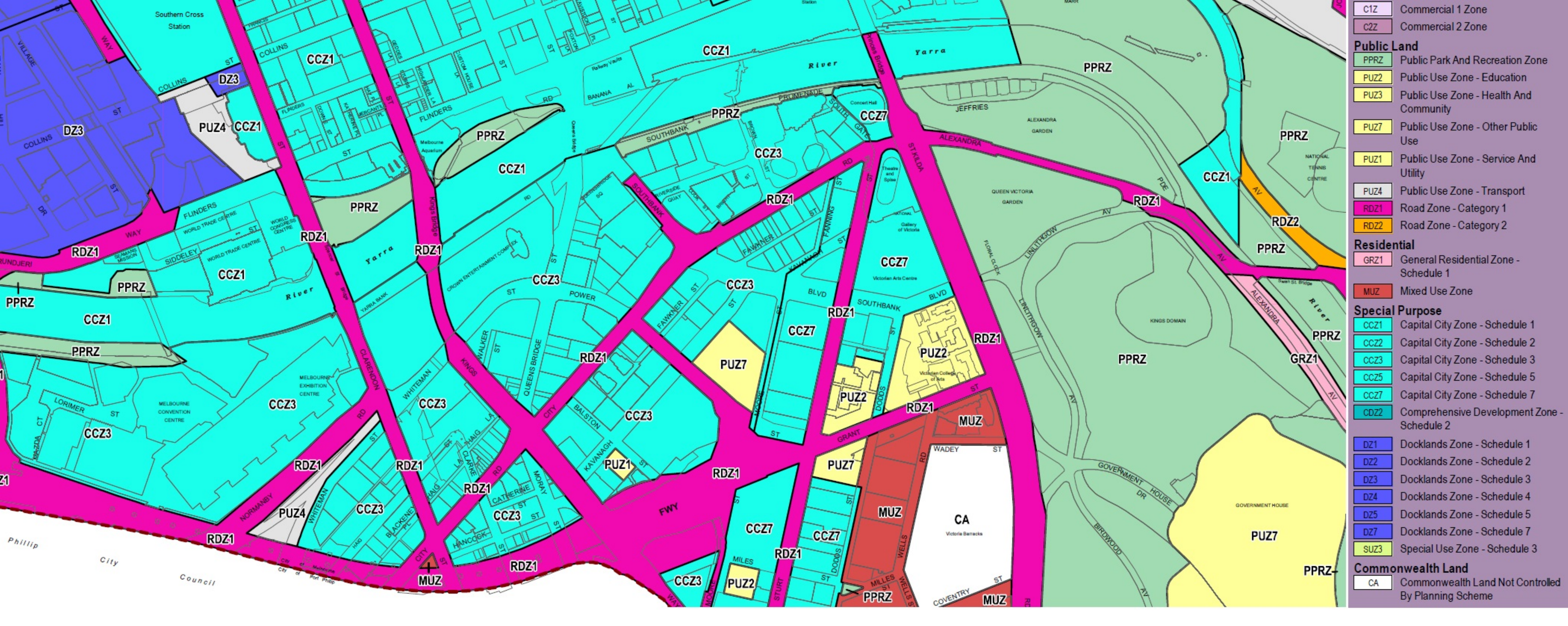}
\caption{Melbourne Zoning Map}
\label{fig:melbsouthbankzoningmap}
\end{figure*}

\begin{figure*}[!h]
\centering
\includegraphics[width=17cm]{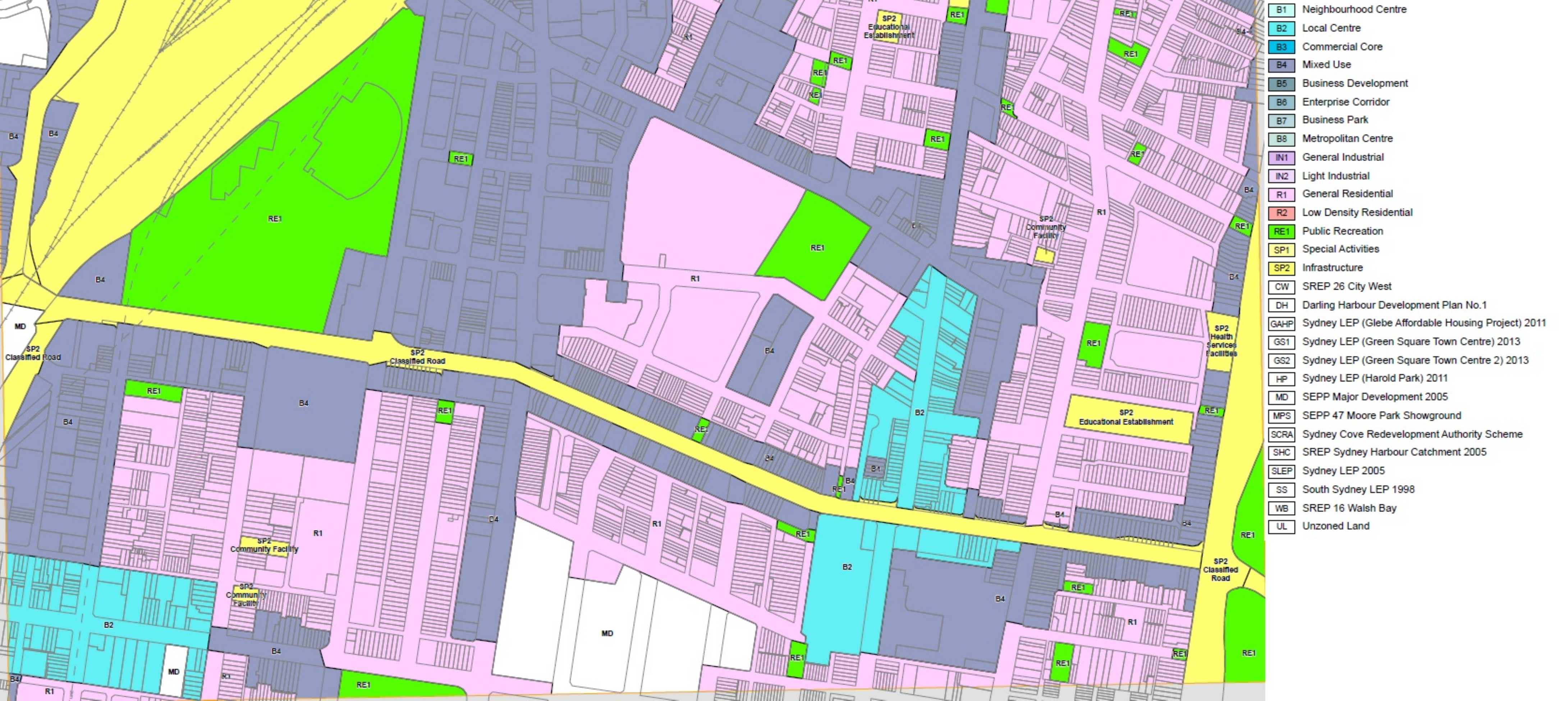}
\caption{Sydney Zoning Map}
\label{fig:sydneyredfernzoningmap}
\end{figure*}

In this section, we validate the extracted clusters of Melbourne and Sydney against actual land use and evaluate the overall methodology. To validate the empirical results obtained so far, the actual land zoning information provided by the city councils of Melbourne and Sydney is utilised. Comparison is made against the type of land zone corresponding to the predicted cluster type. The Melbourne city council divides land into commercial, recreation, residential and mixed-use zones etc. as shown in figure \ref{fig:melbsouthbankzoningmap}. Similarly, the Sydney city council divides land into commercial, mixed-use, residential and industrial zones etc. as shown in figure \ref{fig:sydneyredfernzoningmap}. \newline

In order to determine how closely the identified clusters signify the official land uses, we calculate the degree of overlap between polygon shape of each extracted cluster and the polygon of relevant underlying space in the official zoning map of a city. Table \ref{tab:percentagesofoverlap} shows the percentages of overlap for various types of clusters for both cities. We find that in general there is a high similarity between the official land use types and the predicted land use types from twitter.
\begin{table}[htbp]
\begin{adjustbox}{width=\columnwidth,center}
    \centering
    \begin{tabular}{lcccc}\hline
    
     {\textbf{Predicted Land Use}} & \multicolumn{4}{c}{\textbf{Official Land Use}}\\ 
    
    \multicolumn{1}{c}{} & \multicolumn{1}{c}{\textbf{Business}} & \multicolumn{1}{c}{\textbf{Residential}} & \multicolumn{1}{c}{\textbf{Education}} & \multicolumn{1}{c}{\textbf{Recreation}} \\ \hline
    
    \multicolumn{5}{l}{\textbf{Melbourne}}  \\ 
    
    Business & \multicolumn{1}{c}{\textbf{52.4\%}} & - & - & - \\ 
    
    Residential & - & \multicolumn{1}{c}{\textbf{56\%}} & - & -  \\ 
    
    Education & - & - & \multicolumn{1}{c}{0\%} & - \\ 
    
    Recreation & - & - & - & \multicolumn{1}{c}{\textbf{67.2\%}}   \\ \hline
    
    \multicolumn{5}{l}{\textbf{Sydney}} \\ 
    
    Business & {\textbf{62.5\%}} & - & - & - \\ 
    
    Residential & - & {\textbf{50\%}} & -  & -  \\ 
     
    Education & - & - & {\textbf{53.5\%}} & -  \\ 
    
    Recreation & - & - & - & {\textbf{52.6\%}}  \\ \hline
    \end{tabular}
    \end{adjustbox}
\caption{Overlap between the land use predicted by Twitter activity and the Official land use.}

\label{tab:percentagesofoverlap}

\end{table}
According to Table \ref{tab:percentagesofoverlap}, the business areas in the official zoning maps are identified appropriately by our predicted business clusters with accuracy ranging between 52\% and 62\%. Similarly, the percentage of similarity between the residential land zone and our predicted residential cluster is between 50\% and 56\%. The education cluster of Melbourne is predicted incorrectly and hence its accuracy against actual education zone is 0\%, whereas for Sydney it is 53.5\%. The highest accuracy is observed for recreation zone ranging from 52\% to 67\%. \newline

Not all types of official land zones are necessarily predicted in all cities, as is the case with Melbourne's education cluster. This is because its corresponding extracted cluster (Cluster 3) generated slightly higher error with respect to it. In terms of calculating the overlap between extracted cluster and its underlying actual land zone, the area of underlying zone plays a key part. The cluster itself may have been predicted correctly for its actual use, but its overlap percentage is influenced by the size of the council's land zone area. For instance, the predicted cluster may have a small area owing to concentration of tweets in a limited space but its corresponding actual land zone area may be much larger. This is possible for places like parks where people might gather in a specific part of it and their resulting cluster will be much smaller than the rest of the park, hence yielding smaller overlap percentage. Another issue is the best number of clusters. Current methodology relies on a limited number of clusters and determining the best number of clusters may also be a fruitful area for further research. \newline

Overall, it can be argued that the predicted clusters in both cities predominantly match their respective actual land zones. Thus, the proposed methodology in this study demonstrates that geo-tagged tweets show great promise for land use identification, and this suggests profitable areas for future research in how best to tune this method for maximum accuracy. Although, it has been shown that key official land use types are reasonably represented by our predicted clusters, this whole analysis is still preliminary and needs further investigation. Moreover, it would also be worth validating the empirical results against some survey data when available.

\section{Conclusions}
The availability of massive geo-data generated as a result of ubiquitous use of various online and mobile applications has highlighted its prospective use in various fields, one of which is urban development and planning. In this study, a method has been proposed and demonstrated for identifying and predicting land use using Twitter data from Brisbane and subsequently testing on Melbourne and Sydney. The results indicate that geo-tagged tweets can represent the land use in line with official land zoning to a good extent. The percentage of overlap is found to be as high as 67.2\% for Melbourne's recreation zone followed by 62.5\% for Sydney's business zone. The results of this work are comparable to those of \cite{cranshaw2012livehoods} and \cite{frias2014spectral}, but the cluster generation method used here is simpler yet robust compared to those studies which use a much slower clustering method requiring filtering out of many raw data points. On the other hand, our method can work with virtually any volume of data and without removing any records or losing underlying spatial features.

It should be noted that the precise nature of activities being undertaken by twitter users is not determined as a result of our technique at this stage.
\bibliography{conf.bib}
\end{document}